\documentclass[12pt]{article}

\title{Lower cone distribution functions and set-valued quantiles form Galois connections\footnote{In honour of Y. Kabanov on the occasion of his 70th birthday.}}

\author{
\c{C}a\u{g}{\i}n Ararat\footnote{Bilkent University, Department of Industrial Engineering, Ankara, Turkey, \href{cararat@bilkent.edu.tr}{cararat@bilkent.edu.tr}.} \and Andreas H. Hamel\footnote{Free University of Bozen, Faculty of Economics and Management, Bozen-Bolzano, Italy, \href{mailto:andreas.hamel@unibz.it}{andreas.hamel@unibz.it}.}
}

\date{{\small \today}}

\usepackage{array}
\usepackage{verbatim}
\usepackage{float}
\usepackage{graphicx}
\usepackage{amsmath, amssymb, enumerate, color ,amsthm}
\usepackage[colorlinks=true, linkcolor=blue, citecolor=blue, urlcolor=blue, filecolor=blue]{hyperref}

\newtheorem{theorem}{Theorem}
\newtheorem{corollary}{Corollary}
\newtheorem{remark}{Remark}
\newtheorem{lemma}{Lemma}

\newtheorem{proposition}{Proposition}

\newcommand{\of}[1]{\ensuremath{\left( #1 \right)}}

\newcommand{\abs}[1]{\ensuremath{\left| #1 \right|}}
\newcommand{\cb}[1]{\ensuremath{ \left\{ #1 \right\} }}
\newcommand{\sqb}[1]{\ensuremath{ \left[ #1 \right] }}

\newcommand{\bs}{\backslash}

\newcommand{\pend}{ \hfill $\square$ \medskip}

\newcommand{\E}{\ensuremath{\mathbb{E}}}
\renewcommand{\Pr}{\ensuremath{\mathbb{P}}}
\newcommand{\K}{\ensuremath{\mathcal{K}}}
\renewcommand{\O}{\ensuremath{\Omega}}
\renewcommand{\o}{\ensuremath{\omega}}

\renewcommand{\P}{\ensuremath{\mathcal{P}}}

\newcommand{\A}{\ensuremath{\mathcal{A}}}
\newcommand{\B}{\ensuremath{\mathcal{B}}}

\newcommand{\G}{\ensuremath{\mathcal{G}}}
\newcommand{\X}{\ensuremath{\mathcal{X}}}

\newcommand{\R}{\mathbb{R}}

\newcommand{\N}{\mathrm{I\negthinspace N}}

\newcommand{\cl}{{\rm cl \,}}

\newcommand{\co}{{\rm co \,}}

\newcommand{\triup}{{\rm \vartriangle}}

\newcommand{\Int}{{\rm int\,}}

\definecolor{color0}{gray}{.50}
\definecolor{color1}{rgb}{0,.2,.8}
\definecolor{color2}{rgb}{1,.2,0}
\definecolor{color3}{rgb}{.8,.5,1}

\begin{document}

\maketitle
\begin{abstract} It is shown that the recently introduced lower cone distribution function and the associated set-valued multivariate quantile generate a Galois connection between a complete lattice of closed convex sets and the interval $[0,1]$. This generalizes the corresponding univariate result. It is also shown that an extension of the lower cone distribution function and the set-valued quantile characterize the capacity functional of a random set extension of the original multivariate variable along with its distribution.
\end{abstract}



\section{Introduction}

Several features of set-valued quantiles for multivariate random variables as introduced in \cite{HamelKostner18JMVA} are investigated and extended. In particular, the lower cone distribution function from \cite{HamelKostner18JMVA} is extended to a function on sets, and it is shown that this extension together with the set-valued quantile forms a Galois connection between $([0,1], \leq)$ and a complete lattice of closed convex sets ordered by $\supseteq$. In the univariate case, a similar result is known (see \cite[Remark~3.1]{DoeringDewitt12ArX}), but apparently not too popular under this label. For example, in the recent work \cite{FaugerasRueschendorf17MA}, the property constituting the Galois connection is stated (formula~(2), p.~5), but the Galois connection is neither mentioned, nor exploited.

Our approach turns two downsides of previous proposals for multivariate quantiles into upsides. First, by using the cone distribution function (instead of the joint distribution function even if the cone is $\R^d_+$), an arbitrary vector order can be dealt with, thus `the absence of a natural ordering of Euclidean spaces of dimension greater than one' (\cite[p. 214]{Serfling02SN} where `natural' apparently has to be understood as `total' in order theoretic terms) is turned into a huge potential for applications in statistical and financial analysis where such an order relation very often is present by default, e.g., generated by a solvency cone. Secondly, the fact that an inverse of a monotone function usually `defines only a {\em correspondence}, that is, a multi-valued or set-valued mapping' (\cite[p.~5]{FaugerasRueschendorf17MA}) is exploited by understanding quantiles as functions mapping into complete lattices of sets which carry a rich (order) structure. It is shown that certain lattices, e.g., generated by the closure operators of the respective Galois connection, characterize features of the underlying random vector. 

Moreover, it is shown that the set-valued quantiles characterize the distribution of a random set extension of the original random variable, thus  the three objects ``distribution of the random set," ``(extended) cone distribution function of the random vector," and ``lattice-valued quantile function" carry the same information. This is very much parallel to the univariate case (compare, for instance, \cite[formula~(4), p.~5]{FaugerasRueschendorf17MA}). The ordering cone enters the definition of the lattice of sets in which the random set extension of the original multivariate random variable takes its values.

\section{Set-up}

The framework and notation of \cite{HamelKostner18JMVA} and, when it comes to concepts from set-valued convex analysis, \cite{HamelEtAl15Incoll} are used. In particular, we consider a vector preorder on $\R^d$ which is generated by a nonempty closed convex cone $\emptyset \neq C \subsetneq \R^d$ by means of 
\[
y \leq_C z \quad \Leftrightarrow \quad z - y \in C
\]
for $y, z \in \R^d$; $C$ is neither assumed to have a non-empty interior, nor be pointed, i.e., $C \cap (-C) = \{0\}$ is not assumed. Thus, the cases $C = \{0\}$ and $C = H^+(w) := \{z \in \R^d \mid w^Tz \geq 0\}$ for $w \in \R^d\bs\{0\}$ are not excluded. In the latter case, $C$ is a (homogeneous) halfspace and $\leq_C$ a total preorder. The (positive) dual of the cone $C$ is
\[
C^+ = \cb{w \in C \mid \forall z \in C \colon w^Tz \geq 0}.
\]
The bipolar theorem yields $C = C^{++} := (C^+)^+$ under the given assumptions. The set 
\[
\G(\R^d, C) = \cb{D \subseteq \R^d \mid D = \cl\co(D + C)}
\]
comprises the closed convex subsets of $\R^d$ which are stable under addition of $C$; the sum $A+B$ is understood in the Minkowski sense with $A + \emptyset = \emptyset + A =\emptyset$ for all $A \subseteq \R^d$. The pair $(\G(\R^d, C), \supseteq)$ is an order complete lattice with the following formulas for inf and sup (see, for example, \cite{HamelEtAl15Incoll}) for sets $\mathcal D \subseteq \G(\R^d, C)$:
\[
\inf_{D \in \mathcal D} D = \cl\co\bigcup_{D \in \mathcal D} D, \qquad \sup_{D \in \mathcal D} D = \bigcap_{D \in \mathcal D} D.
\]

\section{Lower cone distribution functions and quantiles}

Let $(\Omega, \A, \Pr)$ be a probability space and $L^0_d$ the space of all equivalence classes of random variables with values in $\R^d$. For $X \in L^0_d$, $w \in \R^d\bs\{0\}$, the function $F_{X, w} \colon \R^d \to [0,1]$ defined by $F_{X, w}(z) = \Pr\{X \in z - H^+(w)\} = \Pr\{w^TX \leq w^Tz\}$ is called the $w$-distribution function of $X$. If $d = 1$, $C = \R_+$ and $w=1$, then $F_{X, w}$ is the usual cumulative distribution function (cdf) of the univariate random variable $X$. The function $F_{X, C} \colon \R^d \to [0,1]$ defined by
\[
F_{X, C}(z) = \inf_{w \in C^+\bs\{0\}}F_{X, w}(z) =  \inf\cb{\Pr\{X \in z - H^+(w)\} \mid w \in C^+}
\]
is called the lower $C$-distribution function associated to $X$.

If $p \in [0,1]$ and $w \in \R^d\bs\{0\}$,  the set $Q^-_{X, w}(p) = \cb{z \in \R^d \mid F_{X, w}(z) \geq p}$ is called the lower $w$-quantile, and the set
\begin{equation}
\label{Qint}
Q^-_{X, C}(p) = \cb{z \in \R^d \mid F_{X, C}(z) \geq p} = \bigcap_{w \in C^+\bs\{0\}}Q^-_{X, w}(p)
\end{equation}
is called the lower $C$-quantile of $X$. Clearly, $Q^-_{X, w}(0) = Q^-_{X, C}(0) = \R^d$ for all $w \in \R^d\bs\{0\}$.  

If $C = \{0\}$ and hence $C^+ = \R^d$, then $F_{X, C}$ is the Tukey depth function and $Q^-_{X, C}(p)$ is the corresponding depth region. In this case, it might happen that (for continuous distributions, for example) $F_{X, C}(z) \in [0, 1/2]$ for all $z \in \R^d$ which also shows that Tukey's depth function in the case $d =1$ is not a generalization of the univariate cdf which requires $C=C^+ = \R_+$ and has values in $[0,1]$, in general.

A few elementary properties are collected in the following result.

\begin{proposition}
\label{PropBasicProperties}
(a) For each $p \in [0,1]$ and $w\in\R^d\bs\{0\}$, the set $Q^-_{X, w}(p)$ is a closed halfspace or empty or $\R^d$.

(b) For each $p \in [0,1]$, the set $Q^-_{X, C}(p)$ is closed, convex and satisfies $Q^-_{X, C}(p) + C \subseteq Q^-_{X, C}(p)$.

(c) One has $Q^-_{X, C}(p_1) \supseteq Q^-_{X, C}(p_2)$ whenever $0 \leq p_1 \leq p_2 \leq 1$.  Moreover,
\begin{equation}
\label{EqLeftContQuantiles}
\forall p \in [0,1] \colon Q^-_{X,C}(p) = \bigcap_{0 \leq q < p}Q^-_{X,C}(q).
\end{equation}

(d) The function $F_{X, C} \colon \R^d \to [0,1]$ is quasiconcave, upper semicontinuous and monotone nondecreasing with respect to $\leq_C$.  
\end{proposition}

{\sc Proof.} (a) This is a consequence of the monotonicity and the upper semicontinuity of $F_{X, w}$ as the cumulative distribution function of the univariate random variable $w^TX$.

(b) This follows from (a) since $Q^-_{X, C}(p)$ is the intersection of the closed halfspaces $Q^-_{X, w}(p)$ for $w \in C^+\bs\{0\}$ (possibly empty or $\R^d$).

(c) The sets $Q^-_{X,C}(p)$ are nested by definition. Hence $\supseteq$ holds in \eqref{EqLeftContQuantiles}. For the contrary, assume $z \in Q^-_{X, w}(p)$. If $z$ is not an element of the right hand side of \eqref{EqLeftContQuantiles}, then there is $0 \leq q < p$ such that $z \not\in Q^-_{X,C}(q)$. This would imply $F_{X, C}(z) < q < p \leq F_{X, C}(z)$, a contradiction.

(d) The first two properties follow since a function with convex and closed upper level sets is quasiconcave and upper semicontinuous while monotonicity is a straightforward consequence of the definitions of $C^+$, $F_{X, w}$, $F_{X, C}$. \pend

Proposition \ref{PropBasicProperties} yields $Q^-_{X,C}(p) \in \G(\R^d, C)$ which means that $Q^-_{X,C}$ can be seen as a function mapping from $[0,1]$ into the complete lattice $(\G(\R^d, C), \supseteq)$. Therefore, \eqref{EqLeftContQuantiles} can be written as $Q^-_{X,C}(p) = \sup_{0 \leq q < p}Q^-_{X,C}(q)$ where the supremum is understood in $(\G(\R^d, C), \supseteq)$. In this sense, $Q^-_{X,C}$ is left-continuous. To summarize, the quantile function $p \mapsto Q^-_{X,C}(p)$ is the non-decreasing, left-continuous $\G(\R^d, C)$-valued inverse of the lower $C$-distribution function $F_{X,C}$ (in the sense of, e.g., \cite[Definition 1]{DrapeauHamelKupper16SVVAN}). This provides a complete analog to the univariate case. The left-continuity of $Q^-_{X,C}$ yields $Q^-_{X,C}(1) = \bigcap_{0 \leq q < 1}Q^-_{X,C}(q)$, and this set can be non-empty. Since $Q^-_{X,C}(0) = \R^d$ is the obvious choice, $Q^-_{X,C}$ is well-defined on $[0,1]$ by \eqref{Qint}. Even for the univariate case, it has been observed that `leaving out the probabilities 0 and 1 is artificial' (\cite[Remark~3.1]{DoeringDewitt12ArX}).

Proposition \ref{PropBasicProperties} (d) can be considered as an extension of \cite[Proposition~1]{RousseeuwRuts99Met} which works for the Tukey depth function (but even for an arbitrary positive measure).

The next result, stated for notational convenience, prepares a continuity result for $F_{X, C}$.

\begin{lemma}\label{limitlemma}
	Let $(a_n)_{n\in\mathbb{N}}, (b_n)_{n\in\mathbb{N}}$ be convergent sequences in $\R$ with limits $a,b\in\R$, respectively. If $a\neq b$, then,
	\[
	\lim_{n\rightarrow\infty}1_{(-\infty,a_n]}(b_n)=1_{(-\infty,a]}(b).
	\]
	\end{lemma}

{\sc Proof.} Suppose that $a<b$ so that $1_{(-\infty,a]}(b)=0$. Let $\varepsilon=\frac{b-a}{3}$. There exists $n_0\in\mathbb{N}$ such that $\abs{a_n-a}\leq \varepsilon$ and $\abs{b_n-b}\leq \varepsilon$ for every $n\geq n_0$. In particular, $a_n\leq a+\varepsilon<b-\varepsilon\leq b_n$ so that $1_{(-\infty,a_n]}(b_n)=0$ for every $n\geq n_0$. Hence, $\lim_{n\rightarrow\infty}1_{(-\infty,a_n]}(b_n)=0=1_{(-\infty,a]}(b).$ The case $a > b$ can be treated by similar arguments.
 \pend

 \begin{remark}
 The condition $a\neq b$ in Lemma~\ref{limitlemma} cannot be omitted. As a counterexample, let $a=b=0$ and $a_n=-\frac1n, b_n=\frac1n$ for every $n\in\mathbb{N}$. Note that $1_{(-\infty,a]}(b)=1$ and $1_{(-\infty,a_n]}(b_n)=0$ for every $n\in\mathbb{N}$. Hence, $\lim_{n\rightarrow\infty}1_{(-\infty,a_n]}(b_n)=0\neq 1=1_{(-\infty,a]}(b)$.
 \end{remark}

\begin{proposition} If the distribution of $w^TX$ under $\Pr$ is such that $\Pr\{w^TX = r\} = 0$ for each $w \in C^+\bs\{0\}$ and each $r \in \R$, then $F_{X,C}$ is continuous. In particular, $F_{X,C}$ is continuous whenever $X$ is a continuous random vector.
\end{proposition}

{\sc Proof.} Let $B=C^+\cap \mathbb{S}^{d-1}$, where $\mathbb{S}^{d-1}$ is the unit sphere in $\R^d$. Note that $B$ is a base for $C^+$ in the sense that every $\tilde{w}\in C^+\setminus\cb{0}$ can be written in the form $\tilde{w}=rw$ for some unique $r>0$ and unique $w\in B$, and we have $F_{X,\tilde{w}}(z)=F_{X,w}(z)$ for every $z\in\R^d$. It follows that
\begin{align}\label{baseinf}
F_{X,C}(z)=\inf_{w\in B}F_{X,w}(z)
\end{align}
for every $z\in \R^d$. Moreover $B$ is a compact set.

By Proposition~\ref{PropBasicProperties} (d), it suffices to show that $F_{X,C}$ is lower semicontinuous. We fix $p\in [0,1)$ and show that the lower level set $L(p)=\cb{z\in\R^d\mid F_{X,C}(z)\leq p}$ is closed. To that end, let $(z_n)_{n\in\mathbb{N}}$ be a convergent sequence in $L(p)$ with limit $\bar{z}\in\R^d$. Let $\varepsilon>0$. By \eqref{baseinf}, for every $n\in\mathbb{N}$, there exists $w_n\in B$ such that $F_{X,w_n}(z_n)< F_{X,C}(z_n)+\varepsilon\leq p+\varepsilon$. Since $(w_n)_{n\in\mathbb{N}}$ is a sequence in the compact set $B$, by Bolzano-Weierstrass theorem, there exists a convergent subsequence $(w_{n_k})_{k\in\mathbb{N}}$ of it, say, with limit $\bar{w}\in B$. Hence, $\lim_{k\rightarrow\infty}w_k^Tz_k=\bar{w}^T\bar{z}$, and applying Lemma~\ref{limitlemma} gives
\[
\lim_{k\rightarrow\infty}1_{(-\infty,w_{n_k}^Tz_{n_k}]}(w_{n_k}^TX(\omega))=1_{(-\infty,\bar{w}^T \bar{z}]}(\bar{w}^TX(\omega))
\]
for every $\omega\in\Omega$ such that $\bar{w}^TX(\omega)\neq \bar{w}^T\bar{z}$. By assumption, one has $\Pr\{\bar{w}^TX=\bar{w}^T\bar{z}\}=0$. Hence,
\[
\lim_{k\rightarrow\infty}1_{(-\infty,w_{n_k}^Tz_{n_k}]}(w_{n_k}^TX)=1_{(-\infty,\bar{w}^T\bar{z}]}(\bar{w}^TX)\text{ almost surely}.
\]
Therefore, by dominated convergence theorem,
\begin{align*}
\lim_{k\rightarrow\infty}F_{X,w_{n_k}}(z_{n_k})&=\lim_{k\rightarrow\infty}\Pr\{w_{n_k}^T X\leq w_{n_k}^T z_{n_k}\}\\
&=\lim_{k\rightarrow\infty}\E\sqb{1_{(-\infty,w_{n_k}^Tz_{n_k}]}(w_{n_k}^TX)}\\
&=\E\sqb{1_{(-\infty,\bar{w}^T\bar{z}]}(\bar{w}^TX)}\\
&=\Pr\{\bar{w}^TX\leq \bar{w}^T\bar{z}\}=F_{X,\bar{w}}(\bar{z}).
\end{align*}
Hence, $F_{X,\bar{w}}(\bar{z})\leq p+\varepsilon$. Since $\varepsilon>0$ is arbitrary, we conclude that $F_{X,C}(\bar{z})\leq F_{X,\bar{w}}(\bar{z})\leq p$. So $\bar{z}\in L(p)$. Hence, $L(p)$ is a closed set.
\pend

There is an alternative way of writing $Q^-_{X,C}$: the result in Theorem~\ref{ThmCQuantileDual} below is inspired by \cite[Propositions~2~\&~6]{RousseeuwRuts99Met}. The proof is prepared by the following lemma which should be known and is implicitly part of the proof of \cite[Theorem~2.11]{ZuoSerfling00AS}.

\begin{lemma}
\label{LemWDistDual}
Let $p\in[0,1]$. For every $w \in \R^d\bs\{0\}$ and every $z \in \R^d$ with $\Pr\{X \in z - H^+(w)\} < p$ there exists $y \in z + \Int H^+(w)$ such that $\Pr\{X \in y - \Int H^+(w)\} < p$.
\end{lemma}

{\sc Proof.} Let us fix $w \in \R^d\bs\{0\}$ and $z \in \R^d$ with $\Pr\{X \in z - H^+(w)\} < p$ and take $\bar z \in \R^d$ with $w^T\bar z = 1$, which exists since $w \neq 0$. Then, $s\bar z \in \Int H^+(w)$ for all $s > 0$. Let us define $y_n = z + \frac{1}{n}\bar z \in z + \Int H^+(w)$ for each $n \in\mathbb{N}$. Then, for every $n\in\mathbb{N}$,
\[
w^Ty_n = w^T\of{z + \frac{1}{n}\bar z} = w^Tz + \frac{1}{n},
\]
so that $w^Ty_{n+1} < w^Ty_n$ and $\lim_{n \to \infty} w^Ty_n = w^Tz$. Since $F_{w^TX}$ is right-continuous, it follows that
\[
\Pr\{X \in z - H^+(w)\} = F_{w^TX}(w^Tz) = \lim_{n \to \infty}F_{w^TX}(w^Ty_n) < p,
\]
so there is $\bar n \in \mathbb{N}$ with
\[
F_{w^TX}(w^Ty_{\bar n}) =  \Pr\{X \in y_{\bar n} - H^+(w)\} < p.
\]
Hence,
\[
\Pr\{X \in y_{\bar n} - \Int H^+(w)\} \leq \Pr\{X \in y_{\bar n} - H^+(w)\} < p,
\]
which proves the claim with $y = y_{\bar n}$. \pend

\begin{theorem}
\label{ThmCQuantileDual}
For all $p \in [0,1]$,
\begin{align*}
Q^-_{X,C}(p) & = \bigcap_{w \in C^+\bs\{0\}}\bigcap_{y \in \R^d}\cb{y + H^+(w) \mid \Pr\cb{X \in y + H^+(w)} > 1-p} \\
	& = \bigcap_{w \in C^+\bs\{0\}}\bigcap_{y \in \R^d}\cb{y + H^+(w) \mid \Pr\cb{X \in y - \Int H^+(w)} < p}.
\end{align*}
\end{theorem}

{\sc Proof.} The two expressions on the right hand side clearly coincide since $\Pr\{X \in y + H^+(w)\} = 1 - \Pr\{X \in y - \Int H^+(w)\}$.

First, assume that $z \not\in \bigcap_{w \in C^+\bs\{0\}}\bigcap_{y \in \R^d}\cb{y + H^+(w) \mid \Pr\cb{X \in y - \Int H^+(w)} < p}$. Then, there exist $w \in C^+\bs\{0\}$, $y \in \R^d$ such that $\Pr\cb{X \in y - \Int H^+(w)} < p$ and $z \not\in y + H^+(w)$. It follows that $z \in y - \Int H^+(w)$ which implies $z - H^+(w) \subseteq y - \Int H^+(w)$, so
\[
\Pr\{X \in z - H^+(w)\} \leq \Pr\{X \in y - \Int H^+(w)\} < p,
\]
hence $z \not\in Q^-_{X,C}\of{p}$. 

Therefore, $Q^-_{X,C}\of{p} \subseteq \bigcap_{w \in C^+\bs\{0\}}\bigcap_{y \in \R^d}\cb{y + H^+(w) \mid \Pr\cb{X \in y - \Int H^+(w)} < p}$.

Conversely, assume that
\[
\bar z \notin Q^-_{X,C}(p) = \bigcap_{w \in C^+\bs\{0\}}Q^-_{X, w}\of{p} = \bigcap_{w \in C^+\bs\{0\}} \cb{z \in \R^d \mid  F_{X,w}(z) \geq p}.
\]
Then, there exists $w \in C^+$ such that $F_{X,w}(\bar z) = \Pr\{X \in \bar z - H^+(w)\} < p$. Lemma \ref{LemWDistDual} yields the existence of $\bar y \in \bar z + \Int H^+(w)$ satisfying $\Pr\{X \in \bar y - \Int H^+(w)\} < p$. If 
\[
\bar z \in \bigcap\limits_{w \in C^+}\bigcap\limits_{y \in \R^d}\cb{y + H^+(w) \mid \Pr\cb{X \in y - \Int H^+(w)} < p}
\]
would be true, then also $\bar z \in \bar y + H^+(w)$ and
\[
\bar z \in \of{\bar y -\Int H^+(w)} \cap \of{\bar y + H^+(w)},
\]
which is a contradiction. So, $\bar z \not\in \bigcap_{w \in C^+}\bigcap_{y \in \R^d}\cb{y + H^+(w) \mid  \Pr\cb{X \in y - \Int H^+(w)} < p}$. This shows $Q^-_{X,C}\of{p} \supseteq \bigcap_{w \in C^+}\bigcap_{y \in \R^d}\cb{y + H^+(w) \mid \Pr\cb{X \in y - \Int H^+(w)} < p}$. 
\pend

\begin{remark}
Assume that there is $\bar z \in C$ such that $w^T \bar z > 0$ for all $w \in C^+\bs\{0\}$. Then, the set $B^+ = \cb{w \in C^+ \mid w^T \bar z = 1} \subseteq C^+$ is a base of $C^+$ (every $w \in C^+\bs\{0\}$ can be represented uniquely as $w = sb$ with $b \in B^+$, $s >0$). If this is the case, then intersections such as in formula \eqref{Qint} and the one in Theorem \ref{ThmCQuantileDual} need to run only over $B^+$ instead of $C^+$ since $H^+(sw) = H^+(w)$ for all $w \in C^+$ and $s >0$.

If $d=1$, $C = \R_+$, then $B^+ = \{1\}$ is such a base, and the formula in Theorem \ref{ThmCQuantileDual} breaks down to 
\[
Q^-_{X,C}(p) = \sup\{r \in \R \mid \Pr\{X < r\} < p\} + \R_+
\]
while \eqref{Qint} becomes $Q^-_{X,C}(p) = \inf\{s \in \R \mid \Pr\{X \leq s\} \geq p\} + \R_+$ which re-produces well-known formulas for the univariate lower quantile, see \cite[p.~207]{FoellmerSchied11Book}. 
\end{remark}

\section{Galois connections}

Let $\P(\R^d)$ denote the power set of $\R^d$ (including $\emptyset$) and $\psi\colon\R^d \to \R\cup\{\pm\infty\}$ be an extended real-valued function. The function $\psi^\triup \colon \P(\R^d) \to \R\cup\{\pm\infty\}$ defined by
\[
\psi^\triup(D) = \inf_{z \in D} \psi(z)
\]
is called the inf-extension of $\psi$, where $\psi^\triup(\emptyset) = +\infty$ is understood. 

First, we apply this concept to $\psi = F_{X,C}$. The following result collects a few properties of $F_{X,C}^\triup$ which are basically inherited from $F_{X, C}$.

\begin{proposition}
\label{PropInfFProps}
(a) $X\mapsto F^\triup_{X, C}(D)$ is monotone: $X^1 \leq_{C} X^2$ almost surely implies $F^\triup_{X^2, C}(D) \leq F^\triup_{X^1, C}(D)$ for every $D \in \P(\R^d)$.

(b) $D\mapsto F^\triup_{X, C}(D)$ is monotone: $D_1 \supseteq D_2$ implies $F^\triup_{X, C}(D_1) \leq F^\triup_{X, C}(D_2)$ for every $X \in L^0_d$.

(c) $F^\triup_{X, C}(\{z\} + C) = F^\triup_{X, C}(\{z\}) =F_{X, C}(z)$ for every $z \in \R^d$.
\end{proposition}

{\sc Proof.} (a) follows from the monotonicity of $X \mapsto F_{X, C}(z)$. (b) follows by the construction of $F^\triup_{X, C}$. (c) follows from the monotonicity of $z \mapsto F_{X, C}(z)$. \pend

The following proposition prepares a new feature.

\begin{proposition}
\label{PropPreInfStable}
For every $D \in \P(\R^d)$,
\[
F^\triup_{X, C}(\cl\co(D + C)) = F^\triup_{X, C}(D). 
\]
\end{proposition}

{\sc Proof.} By the definition of $F^\triup_{X, C}$, $\leq$ is certainly true. Since $F_{X, C}$ is monotone with respect to $\leq_C$ (see Proposition~\ref{PropBasicProperties} (d)),
\[
\forall z \in C, \; \forall x \in D \colon F_{X, C}(x) \leq F_{X, C}(x + z),
\]
hence $F^\triup_{X, C}(\{x\}) \leq F^\triup_{X, C}(\{x\}+C)$ and therefore, $F^\triup_{X, C}(D) \leq F^\triup_{X, C}(D + C)$. This gives $F^\triup_{X, C}(D) = F^\triup_{X, C}(D + C)$.

Next, take $z^1, \ldots, z^m \in D$ and $s_1,\ldots,s_m \in [0,1]$ with $m\in\mathbb{N}$ such that $\sum_{i=1}^m s_i =1$. Set $z = \sum_{i=1}^m s_iz^i$. The quasiconcavity of $F_{X, C}$ yields
\[
F_{X, C}(z) \geq \min\{F_{X, C}(z^1), \ldots, F_{X, C}(z^m)\} \geq F^\triup_{X, C}(D).
\]
This proves $F^\triup_{X, C}(\co D) = F^\triup_{X, C}(D)$.

Finally, take a sequence $(z^n)_{n\in\mathbb{N}}$ in $D$ which converges to some $z \in \R^d$. Then, the upper semicontinuity of $F_{X, C}$ produces
\[
F_{X, C}(z) \geq \limsup_{n\rightarrow\infty} F_{X, C}(z^n) \geq \limsup_{n\rightarrow\infty} F^\triup_{X, C}(D) = F^\triup_{X, C}(D), 
\]
hence $F^\triup_{X, C}(\cl D) \geq F^\triup_{X, C}(D)$ which gives ``$=$." \pend

As a result of Proposition~\ref{PropPreInfStable}, it is enough to consider the inf-extension $F^\triup_{X,C}$ as a function on $\G(\R^d,C)$ rather than the whole power set $\P(\R^d)$, which is done in the sequel.

\begin{corollary}
\label{CorInfStable}
(a) For every $\mathcal D \subseteq  \G(\R^d, C)$, one has
\begin{equation}
\label{EqInfStability}
F^\triup_{X, C}\of{\inf_{D \in \mathcal D} D} = \inf_{D \in \mathcal D} F^\triup_{X, C}(D).
\end{equation}
(b) For $A \subseteq \R^d$ and $w \in C^+$, one has
\begin{equation}
\label{EqInfStabilityW}
\inf_{y \in A} F_{w^TX}(w^T y) = F_{w^TX}\of{\inf_{y \in A} w^Ty}.
\end{equation}
\end{corollary}

{\sc Proof.} (a) Since $D \subseteq \inf_{D \in \mathcal D} D = \cl\co\bigcup_{D \in \mathcal D} D$ for all $D \in \mathcal D$, ``$\leq$" certainly is true. Take $z \in \bigcup_{D \in \mathcal D} D$. Then, there is $D' \in \mathcal D$ with $z \in D'$, hence $F_{X, C}(z) \geq  \inf_{D \in \mathcal D} F^\triup_{X, C}(D)$ which in turn implies  $F^\triup_{X, C}(\bigcup_{D \in \mathcal D}D) \geq \inf_{D \in \mathcal D} F^\triup_{X, C}(D)$. Proposition~\ref{PropPreInfStable} now produces $F^\triup_{X, C}(\inf_{D \in \mathcal D} D) = \inf_{D \in \mathcal D} F^\triup_{X, C}(D)$. 

(b) Setting $C = H^+(w)$, $D(y) = y + H^+(w)$ for $y \in A$, $\mathcal D = \cb{D(y)\mid y \in A}$ and observing $F^\triup_{X, H^+(w)}(A) = \inf_{y \in A}F_{w^TX}(w^T y)$ as well as
\[
F^\triup_{X, w}\of{\inf_{D \in \mathcal D} D} = F^\triup_{X, w}\of{\cl\bigcup\limits_{y \in A} (y + H^+(w))} = F^\triup_{X, w}(A) = \inf_{y \in A} F_{w^TX}(w^T y)
\]
(see Proposition \ref{PropPreInfStable}) one gets \eqref{EqInfStabilityW} as a special case of \eqref{EqInfStability}.
\pend

\medskip Equation \eqref{EqInfStability} means that $F^\triup_{X, C}$ preserves infima (meets) as a function from $\G(\R^d, C)$ to $[0,1]$ (see also \cite[7.31 \& 2.26]{DaveyPriestley02SecBook}). This property has been called ``inf-stability" in \cite{CrespiEtAl18ArX}.

\begin{proposition}
\label{PropGalois} (a) For every $D \in \G(\R^d, C)$ and $p \in [0, 1]$,
\[
Q^-_{X,C}(p) \supseteq D \quad \Leftrightarrow \quad p \leq F^\triup_{X, C}(D).
\]

(b) The two compositions $F^\triup_{X, C} \circ Q^-_{X,C} \colon [0,1] \to \G(\R^d, C)$ and $Q^-_{X,C} \circ F^\triup_{X, C} \colon \G(\R^d, C) \to [0,1]$ are closure operators (extensive, increasing and idempotent).

(c) The set
\[
\mathcal F_X(\R^d, C) = \cb{D \in \G(\R^d, C) \mid D = (Q^-_{X,C} \circ F^\triup_{X, C})(D)}
\]
is a complete lattice with respect to $\supseteq$.

(d) One has
\begin{align}
\forall p \in [0,1] \colon \quad Q^-_{X,C}(p) & = \inf\cb{D \in \G(\R^d, C) \mid F^\triup_{X, C}(D) \geq p} \label{FQdet}\\
\forall D \in \G(\R^d, C) \colon  \quad  F^\triup_{X, C}(D)&  = \sup\cb{p \in [0,1] \mid D \subseteq Q^-_{X,C}(p)}.\label{QFdet}
\end{align}
\end{proposition}

{\sc Proof.} (a) is straightforward and can be checked by using the definitions of $F^\triup_{X, C}$ and $Q^-_{X,C}$. (b) follows from the theory of Galois connections (see  \cite[Chapter~7]{DaveyPriestley02SecBook}). (c) follows from the Knaster-Tarski theorem since $(\G(\R^d, C), \supseteq)$ is a complete lattice and $\mathcal F_X(\R^d, C)$ is the set of fixed points of the composition $Q^-_{X,C} \circ F^\triup_{X, C} \colon \G(\R^d, C) \to \G(\R^d, C)$. (d) also follows from the theory of Galois connections \cite{DaveyPriestley02SecBook}. \pend

\medskip Proposition \ref{PropGalois} (a) establishes the fact that $F^\triup_{X, C}$ and $Q^-_{X,C}$ form a Galois connection between the two complete lattices $(\G(\R^d, C), \supseteq)$ and $([0,1], \leq)$ where $F_{X,C}^\triup$ is the upper adjoint and $Q_{X,C}^-$ the lower adjoint. This means that $F^\triup_{X, C}$ and $Q^-_{X,C}$ determine each other; they carry the same information.

\begin{remark}\label{totalrank}
The complete lattice $\mathcal F_X(\R^d, C)$ can be generated in a different, but related way. Using the notation of \cite{CrespiEtAl18ArX}, we set $\Psi=\cb{F_{X,C}}$ and define the $\Psi$-closure of $D\in\G(\R^d,C)$ by
\[
\cl_{\Psi}(D)=\cb{z\in\R^d\mid F_{X,C}(z)\geq F_{X,C}^{\triup}(D)}.
\]
Now, one has
\[
\cl_{\Psi}= Q^-_{X,C}\circ F^\triup_{X, C}.
\]
Therefore, the set of all fixed points of the closure operator $Q^-_{X,C} \circ F^\triup_{X, C}$ coincides with the complete lattice generated by the singleton $\Psi=\{F_{X,C}\}$ via $D = \cl_{\psi}(D)$. The relation defined by
\[
z^1 \leq_\Psi z^2 \quad \Leftrightarrow \quad  F_{X,C}(z^1) \leq F_{X,C}(z^2),
\]
is a total order which is extended to $\P(\R^d)$ by
\[
D^1 \preceq_\Psi D^2 \quad \Leftrightarrow \quad  F^\triup_{X,C}(D^1) \leq F^\triup_{X,C}(D^2).
\]
\end{remark}

The function $F_{X,C}$ can be understood as a ranking function for multivariate data points in $\R^d$. Such ranking functions are used in statistics, e.g., for outlier detection (see \cite{Serfling10JNS}), and also for decision making (see \cite{Kostner19MMOR}). The function $F_{X,C}^{\triup}$ gives a corresponding ranking for subsets of $\R^d$. 

A different (non-total) order relation can be constructed using the $w$-distribution functions $F_{X,w}$ with $w\in C^+\bs\{0\}$ instead of the lower $C$-distribution function $F_{X,C}$. We consider the family $\Phi=\cb{F_{X,w}\mid w\in C^+\bs\{0\}}$ which induces the relation
\[
z^1 \leq_\Phi z^2 \quad \Leftrightarrow \quad \forall w\in C^+\bs\{0\}\colon F_{X,w}(z^1) \leq F_{X,w}(z^2)
\]
on $\R^d$ which is non-total in general, as well as the set relation
\[
D^1 \preceq_\Phi D^2 \quad \Leftrightarrow\quad  \forall w\in C^+\bs\{0\}\colon  F^\triup_{X,w}(D^1) \leq F^\triup_{X,w}(D^2)
\]
on $\G(\R^d,C)$. Since the infimum over $w \in C^+$ can be taken first on the left hand side of the two scalar inequalities above and then on the right hand side, the relations $\leq_\Psi$ and $\preceq_\Psi$ above turn out to be extensions of $\leq_\Phi$ and $\preceq_\Phi$, respectively: $z^1 \leq_\Phi z^2$ implies $z^1 \leq_\Psi z^2$ and $D^1 \preceq_\Phi D^2$ implies $D^1 \preceq_\Psi D^2$.

Define the $\Phi$-closure of $D \in \G(\R^d,C)$ by
\[
\cl_\Phi(D)=\bigcap_{w\in C^+\bs\{0\}}\cb{z\in\R^d\mid F_{X,w}(z)\geq F^\triup_{X,w}(D)}.
\]
It follows from \cite[Proposition~2.2]{CrespiEtAl18ArX} that $\cl_{\Phi}$ is a closure operator which generates the complete lattice
\[
\mathcal C_X(\R^d,C)=\cb{D\in\G(\R^d,C)\mid D=\cl_\Phi(D)}
\]
with the relation $\supseteq$. The next result characterizes the case where $\cl_\Phi$ coincides with the identity operator, that is, $\mathcal C_X(\R^d,C)=\G(\R^d,C)$.

\begin{theorem}
	\label{PhiIdentity}
	The following are equivalent:
	
	(a) For every $w\in C^+\bs\{0\}$, the cumulative distribution function $F_{w^T X}$ is strictly increasing.
	
	(b) For every $D\in\G(\R^d,C)$, $D=\cl_\Phi(D)$.
\end{theorem}

{\sc Proof.} Suppose that (a) holds and let $D\in\G(\R^d,C)$. Clearly, $D\subseteq\cl_\Phi(D)$. To show the reverse inclusion, let $z\in \cl_\Phi(D)$. Fix $w\in C^+\bs\{0\}$. We have $F_{X,w}(z)\geq F_{X,w}^\triup(D)$. By \eqref{EqInfStabilityW} one has
\[
F_{w^T X}(w^T z)=F_{X,w}(z)\geq F_{X,w}^\triup(D)=\inf_{y\in D}F_{w^T X}(w^T y)=F_{w^T X}\of{\inf_{y \in D}w^T y},
\]
where $F_{w^T X}(-\infty)=0$ is understood. Hence, the strict monotonicity of $F_{w^\mathsf{T}X}$ implies
\[
w^Tz \geq \inf_{y \in D}w^Ty.
\]
Since this holds for every $w\in C^+\bs\{0\}$, one may conlude that $z\in D$. Hence, $\cl_\Phi(D)\subseteq D$.

Conversely, suppose that (b) holds. To get a contradiction, assume that there exists $\bar{w}\in C^+\bs\{0\}$ such that $F_{\bar{w}^T X}$ is not strictly increasing. Hence, there exist $a,b\in\R$ such that $a<b$ and $F_{\bar{w}^T X}(a)=F_{\bar{w}^T X}(b)$. It is clear that one can find $z^a,z^b\in\R^d$ with $\bar{w}^T z^a=a$ and $\bar{w}^T z^b=b$. Let us define
\[
D=\{z^b\}+H^+(\bar{w}).
\]
Clearly, $D\in\G(\R^d,C)$ and $z^b\in D$. We claim that $z^a\in \cl_\Phi(D)$. First, note that
\[
\inf_{y \in D}w^T y =\begin{cases}s\bar{w}^Tz^b=sb & \text{if } w=s\bar{w} \text{ for some } s>0,\\ -\infty &\text{else},\end{cases}
\]
for each $w\in C^+\bs\{0\}$. Hence, if $w = s\bar{w}$ for $s>0$, then
\begin{align*}
F^\triup_{X,w}(D)&=\inf_{y \in D}F_{w^TX}(w^Ty)=\inf_{y \in D}F_{\bar{w}^TX}(\bar{w}^Ty)\\
& = F_{\bar{w}^T X}\of{\inf_{y \in D}\bar{w}^Ty} = F_{\bar{w}^TX}(b) = F_{\bar{w}^TX}(a)=F_{\bar{w}^TX}(\bar{w}^Tz^a)=F_{X,w}(z^a).
\end{align*}
On the other hand, if $w \in C^+\bs\{0\}$ with $w\neq s\bar{w}$ for every $s>0$, then
\[
F^\triup_{X,w}(D)=F_{w^TX}\of{\inf_{y \in D}w^Ty}=F_{w^T X}(-\infty)=0\leq F_{X,w}(z^a).
\]
Hence, the claim follows. However,
\[
\bar{w}^Tz^a=a<b=\bar{w}^Tz^b=\inf_{y \in D}\bar{w}^Ty,
\]
which shows that $z^a\notin D$. Since $D\neq\cl_\Phi(D)$, we get a contradiction to (b). Hence, (a) holds.
\pend

Note that the condition (a) in Theorem~\ref{PhiIdentity} requires, for each $w\in C^+\bs\{0\}$, the continuous part of $F_{w^TX}$ to be strictly increasing although $F_{w^TX}$ may have jumps.

\begin{remark}
It is easy to check that $D\subseteq\cl_{\Phi}(D)\subseteq\cl_\Psi(D)$ for each $D\in\G(\R^d,C)$. Under the conditions of Theorem~\ref{PhiIdentity}, we have $D=\cl_{\Phi}(D)$. In general, $\cl_{\Phi}(D)$ may be a (much) smaller set than $\cl_{\Psi}(D)$. However, if $C=H^+(w)$ for some $w\in\R^d\bs\{0\}$, then $C^+$ is the ray generated by $w$ so that $\cl_{\Phi}(D)=\cl_{\Psi}(D)$ for every $D\in\G(\R^d,C)$. In this case, $\preceq_\Phi$ is also a total order and coincindes with $\preceq_\Psi$.
\end{remark}

\section{The simulation result}

The main question in this section is how the quantile function characterizes the distribution. In the univariate case, one can show that the quantile, taken at a random variable uniformly distributed over $[0,1]$, produces a random variable which has the cumulative distribution function that defines the quantile (compare, for example, \cite[Lemma A.19]{FoellmerSchied11Book}, the ``simulation lemma"). In our setting, quantiles are sets, so plugging in a random variable with values in $[0,1]$ produces a random set. 

Let $U\colon\O\to[0,1]$ be a standard uniform random variable and define a function $\X\colon\O\to\G(\R^d,C)$ by
\[
\X(\o) = Q^-_{X,C}(U(\o))
\]
for every $\omega\in\O$. To be able to talk about the distribution of $\X$ under $\Pr$, we first view $\X$ as a measurable function by equipping $\G(\R^d,C)$ with the $\sigma$-algebra constructed below.

For each $K\subseteq \R^d$, let
\[
\mathcal{D}^K=\cb{D\in\G(\R^d,C)\mid D\cap K=\emptyset}
\]
and
\[
\mathcal{D}_K=\cb{D\in\G(\R^d,C)\mid D\cap K\neq\emptyset},
\]
which is the complement of $\mathcal{D}^K$ in $\G(\R^d,C)$. Let us denote by $\K$ the set of all compact subsets of $\R^d$. Note that the collection $\cb{\mathcal{D}^K\mid K\in\K}$ is a $\pi$-system on $\G(\R^d,C)$ since $\mathcal{D}^{K_1}\cap \mathcal{D}^{K_2}=\mathcal{D}^{K_1\cup K_2}$ and $K_1\cup K_2\in \mathcal{K}$ for every $K_1,K_2\in\mathcal{K}$. Let $\B(\G(\R^d,C))$ be the $\sigma$-algebra generated by $\cb{\mathcal{D}^K\mid K\in\K}$, called the Borel $\sigma$-algebra on $\G(\R^d,C)$; the reader is referred to \cite[Section~1.1]{Molchanov} for a detailed discussion. Clearly, $\B(\G(\R^d,C))$ is also generated by $\cb{\mathcal{D}_K\mid K\in\K}$.

We shall establish the measurability of $\X$ with respect to $\B(\G(\R^d,C))$.

\begin{lemma}\label{measurable}
	The function $\X\colon\Omega\to\G(\R^d,C)$ is measurable with respect to $\A$ and $\B(\G(\R^d,C))$.
	\end{lemma}

{\sc Proof.} Let $K\in \K$. Note that
\begin{align*}
\cb{\X\in \mathcal{D}_K}&=
\cb{\o\in\O\mid Q^-_{X,C}(U(\o))\cap K\neq \emptyset}\\
&=\cb{\o\in\O\mid \cb{z\in\R^d\mid F_{X,C}(z)\geq U(\o)}\cap K\neq \emptyset}\\
&=\cb{\o\in\O\mid \exists z\in K\colon F_{X,C}(z)\geq U(\o)}\\
&=\cb{\o\in\O\mid \max_{z\in K}F_{X,C}(z)\geq U(\o)}\\
&=\cb{U\leq \max_{z\in K}F_{X,C}(z)}\in \A,
\end{align*}
where the fourth equality and the well-definedness of maximum is by the upper semicontinuity of $F_{X,C}$ in Proposition~\ref{PropBasicProperties}(d), and the last equality follows since $U$ is measurable with respect to the Borel $\sigma$-algebra on $[0,1]$ and $\A$. Hence, by \cite[Proposition~I.2.3]{Cinlar}, it follows that $\cb{\X\in \mathcal{D}}$ for every $\mathcal{D}\in\B(\G(\R^d,C))$, that is, $\X$ is measurable.\pend

Thanks to Lemma~\ref{measurable}, $\X$ is a random variable taking values in $\G(\R^d,C)$. Hence, its distribution under $\Pr$ is the probability measure $\Pr\circ \X^{-1}$ on $(\G(\R^d,C),\B(\G(\R^d,C)))$ defined by
\[
\Pr\circ\X^{-1}(\mathcal{D})=\Pr\cb{\X\in \mathcal{D}}
\]
for every $D\in\B(\G(\R^d,C))$. Since $\cb{\mathcal{D}^K\mid K\in\mathcal{K}}$ is a $\pi$-system which generates the $\sigma$-algebra $\B(\G(\R^d,C))$, the distribution of $\X$ is determined by its values on this $\pi$-system; see \cite[Proposition~I.3.7]{Cinlar}, for instance. Since $\Pr\cb{\X\in \mathcal{D}^K}=1-\Pr\cb{\X\in\mathcal{D}_K}$ for every $K\in\K$, the distribution of $\X$ is also determined by the so-called \emph{capacity functional} $T_{\X}\colon\mathcal{K}\to [0,1]$ defined by
\[
T_\X (K) = \Pr\cb{\X\in \mathcal{D}_K}=\Pr\cb{\X\cap K\neq \emptyset}
\]
for each $K\in\K$.

\begin{proposition}
\label{simthm}
The lower $C$-distribution function $F_{X,C}\colon\R^d\to[0,1]$ and the capacity functional $T_\X$ of the set-valued random variable $\X$ determine each other.
\end{proposition}

{\sc Proof.} Let $K\in\K$. Following the calculation in the proof of Lemma~\ref{measurable}, we have
\[
T_{\X}(K)=\Pr\cb{U\leq \max_{z\in K}F_{X,C}(z)}=\max_{z\in K}F_{X,C}(z)
\]
since $U$ has the standard uniform distribution. Hence, $F_{X,C}$ determines $T_\X$.

Conversely, let $z\in\R^d$. The above calculation yields $F_{X,C}(z)=T_\X(\cb{z})$. Hence, $T_\X$ determines $F_{X,C}$.
\pend

\medskip Proposition~\ref{simthm} together with \eqref{FQdet}, \eqref{QFdet} implies that the lower $C$-quantile $Q^-_{X,C}$, the lower $C$-distribution function $F_{X,C}$, the inf-extension $F^{\triup}_{X,C}$, the capacity functional $T_\X$, and the distribution $\Pr\circ\X^{-1}$ determine each other.

\end{document}